\documentclass[12pt]{article}
\usepackage{exscale,relsize}
\usepackage{lscape}
\usepackage{amsmath}
\usepackage[sort,nocompress]{cite}%increasing order for citations
\usepackage{amsfonts}
\usepackage[hidelinks]{hyperref}
\usepackage{amssymb}
\usepackage{calc}
\usepackage{theorem}
\usepackage{pifont}      %needed by dingautolist
\usepackage{array}
\usepackage{mathtools} % for ceiling function
\DeclarePairedDelimiter{\ceil}{\lceil}{\rceil}

\usepackage{color}

\usepackage{graphicx} %Enable figure
\usepackage{float} % Figure positioning
\usepackage{subcaption} % This package includes subfigure command

\newcommand{\N}{\mathbb{N}}   % natrual numbers

\newcommand{\dist}{\ensuremath{\operatorname{dist}}} %distance function

\newcommand{\norm}[1]{\left\lVert#1\right\rVert} %norm
\newcommand{\ip}[2]{\langle#1,#2\rangle} % inner product

\newcommand{\gph}{\ensuremath{\operatorname{gph}}}

\newcommand{\R}{\ensuremath{\mathbb R}}
\newcommand{\Rn}{\ensuremath{\mathbb R^n}}

\newcommand{\dom}{\ensuremath{\operatorname{dom}}}

\newcommand{\Id}{\ensuremath{\operatorname{Id}}}

\providecommand{\BB}[2]{\mathbb{B}(#1;#2)}

\newcommand{\prox}{\ensuremath{\operatorname{Prox}}}
\newcommand{\proj}{\ensuremath{\operatorname{Proj}}}

\newtheorem{theorem}{Theorem}[section]
\newtheorem{lemma}[theorem]{Lemma}
\newtheorem{fact}[theorem]{Fact}
\newtheorem{corollary}[theorem]{Corollary}
\newtheorem{proposition}[theorem]{Proposition}
\newtheorem{defn}[theorem]{Definition}

\theoremstyle{plain}{\theorembodyfont{\rmfamily}
}
\theoremstyle{plain}{\theorembodyfont{\rmfamily}
}
\theoremstyle{plain}{\theorembodyfont{\rmfamily}
}
\theoremstyle{plain}{\theorembodyfont{\rmfamily}
\newtheorem{example}[theorem]{Example}}
\theoremstyle{plain}{\theorembodyfont{\rmfamily}
\newtheorem{remark}[theorem]{Remark}}
\theoremstyle{plain}{\theorembodyfont{\rmfamily}
}
\def\proof{\noindent{\it Proof}. \ignorespaces}
\def\endproof{\ensuremath{\quad \hfill \blacksquare}}

\newcommand{\pluss}{{\hskip1pt \raise1pt\vbox{\hrule width6pt \vskip1pt
\hrule width6pt}\kern-4pt{\lower1pt\hbox{\vrule height6pt \kern1pt\vrule
height6pt}}\hskip5pt}}

\newcommand{\argmin}{\mathop{\rm argmin}\limits}

\begin{document}

\title{Malitsky-Tam forward-reflected-backward splitting method for nonconvex minimization problems}

\author{
         Xianfu\ Wang\thanks{
                 Mathematics, University of British Columbia, Kelowna, B.C.\ V1V~1V7, Canada.
                 E-mail: \href{mailto:shawn.wang@ubc.ca}{\texttt{shawn.wang@ubc.ca}}.} and
         Ziyuan Wang\thanks{
                 Mathematics, University of British Columbia, Kelowna, B.C.\ V1V~1V7, Canada.
                 E-mail: \href{mailto:ziyuan.wang@alumni.ubc.ca}{\texttt{ziyuan.wang@alumni.ubc.ca}}.}
                 }

\date{November 16, 2021}

\maketitle
\begin{abstract} \noindent We extend the Malitsky-Tam forward-reflected-backward~(FRB) splitting method for inclusion problems of monotone operators to nonconvex minimization problems. By assuming the
generalized concave Kurdyka-\L ojasiewicz~(KL) property of a quadratic regularization of the objective, we show that the FRB method converges globally to a stationary point of the objective and enjoys finite length property. The sharpness of our approach is guaranteed by virtue of the exact modulus associated with the generalized concave KL property. Numerical experiments suggest that FRB is competitive compared to the Douglas-Rachford method and the Bo{\c{t}}-Csetnek inertial Tseng's method.
\end{abstract}

\noindent {\bfseries 2010 Mathematics Subject Classification:}
Primary 49J52, 90C26; Secondary 26D10.

\noindent {\bfseries Keywords:} Generalized concave Kurdyka-\L ojasiewicz property, proximal mapping, Malitsky-Tam forward-reflected-backward splitting method, merit function, global convergence, nonconvex optimization.

	\section{Introduction}

%In the full convex setting where $f$ and $g$ are both convex, there has been an intensive study on splitting methods % for solving~(\ref{optimization prblem}) and the more
Consider the general inclusion problem of maximally monotone operators:
\begin{equation}\label{inclusion problem}
	\text{find }x\in\Rn\text{ such that }0\in A(x)+B(x),
\end{equation}
where $A:\Rn\to 2^{\Rn}$ and $B:\Rn\to\Rn$ are (maximally) monotone operators with $B$ being Lipschitz continuous. %which clearly specifies to~(\ref{optimization prblem}) with $A=\partial f$ and $B=\nabla g$.
Successful splitting methods for the inclusion problem~(\ref{inclusion problem}) are
numerous; see, e.g., \cite[Chapters 26, 28]{BC}, including the forward-backward method~\cite{combettes2005signal} and Tseng's method~\cite{tseng2000modified}. Recently, Malitsky and Tam~\cite{malitsky2020FRB} proposed a \emph{forward-reflected-backward} (FRB) splitting method for solving~(\ref{inclusion problem}). Let $\lambda>0$ and $J_{\lambda A}=(\Id+\lambda A)^{-1}$ be the resolvent of operator $\lambda A$; see,~e.g.,~\cite{BC}. Given initial points $x_0,x_{-1}\in\Rn$, the Malitsky-Tam FRB method with a
fixed step-size (see~\cite[Remark 2.1]{malitsky2020FRB}) iterates
\begin{equation}\label{FRB iteration operator form}
x_{k+1}=J_{\lambda A}(x_k-2\lambda B(x_k)+\lambda B(x_{k-1})),
\end{equation}
with only one forward step computation required per iteration compared to Tseng's method.
In~\cite{malitsky2020FRB} Malitsky and Tam showed that their
FRB given by \eqref{FRB iteration operator form} converges to a solution of~(\ref{inclusion problem}).

%Clearly, FRB serves as a favorable candidate for solving~(\ref{optimization prblem}) in the full convex setting, especially when the forward step computation is expensive.

\emph{In this paper, we extend the Malitsky-Tam FRB splitting method to solve
%the nonconvex setting.
the following structured possibly nonconvex minimization problem:
\begin{equation}\label{optimization prblem}
	\min_{x\in\Rn} F(x)=f(x)+g(x),		
\end{equation}
where $f:\Rn\to\overline{\R}=(-\infty,\infty]$ is proper lower semicontinuous (lsc),
and $g:\Rn\to\R$ has a Lipschitz continuous gradient with constant $L>0$.
This class of optimization problems is of significantly broad interest in both theory and
practical applications, see, e.g.,
\cite{beck2009fast, li2016douglas}.
}

Because both $\partial f, \partial g$ are not necessarily monotone, the techniques developed by Malitsky and Tam
do not apply.
To accomplish our goal, we employ the
generalized concave KL property~(see Definition~\ref{Def: g-concave KL})
recently proposed in~\cite{wang2020} as a key assumption in our analysis, under which the FRB method demonstrates pleasant behavior. We show that the sequence produced by FRB converges globally to a stationary point of~(\ref{optimization prblem}) and has finite length property. The sharpest upper bound on the total length of iterates is described by using the exact modulus associated with the generalized concave KL property. A convergence rates analysis is also carried out. Our approach deviates from the convergence mechanism developed in the original FRB paper~\cite{malitsky2020FRB}, but follows a similar pattern as in many work that devotes to obtaining nonconvex extensions of splitting algorithms; see, e.g.,~\cite{bolte2018first,Attouch2010,boct2016inertial,attouch2013convergence,li2016douglas,li2017peaceman,Banert2019}. Numerical simulations indicate that FRB is competitive compared to the Douglas-Rachford method with a fixed step-size~\cite{li2016douglas} and the Bo{\c{t}}-Csetnek inertial Tseng's method~\cite{boct2016inertial}.

The paper is structured as follows. We introduce notations and basic concepts in Section~\ref{sec:preliminaries}. The FRB convergence analysis is presented in Section~\ref{section: Forward-reflected-backward splitting method for nonconvex problem}. Numerical experiments are implemented in Section~\ref{num}. We end this paper by providing some directions for future research in Section~\ref{section: conclusion and future work}.  	
\section{Notations and preliminaries}\label{sec:preliminaries}

Throughout this paper, $\Rn$ is the standard Euclidean space equipped with inner product~$\ip{x}{y}=x^Ty$ and the Euclidean norm~$\norm{x}=\sqrt{\ip{x}{x}}$ for $x, y\in\Rn$.
Let~$\N=\{0,1,2,3,\ldots\}$ and~$\N^*=\{-1\}\cup\N$. The open ball centered at $\bar{x}$ with radius $r$ is denoted by~$\BB{\bar{x}}{r}$.
The distance function of a subset $K\subseteq \Rn$ is $\dist(\cdot,K):\Rn\rightarrow\overline{\R}=(-\infty,\infty]$,
\[x\mapsto\dist(x,K)=\inf\{\norm{x-y}:y\in K\},\]
where $\dist(x,K)\equiv\infty$ if $K=\emptyset$. For $f:\Rn\to\overline{\R}$ and $r_1,r_2\in[-\infty,\infty]$, we set $[r_1<f<r_2]=\{x\in\Rn:r_1<f(x)<r_2\}$. We say a function $f:\Rn\to\overline{\R}$ is coercive if $\displaystyle\lim_{\norm{x}\to\infty}f(x)=\infty$.

We will use frequently the following subgradients in the nonconvex setting; see, e.g., \cite{rockwets, mor2005variational}.

\begin{defn}\label{Defn:limiting subdifferential}
		Let $f:\Rn\rightarrow\overline{\R}$ be a proper function. We say that
		\begin{itemize}
			\item[(i)] $v\in\Rn$ is a \textit{Fr\'echet subgradient} of $f$ at $\bar{x}\in\dom f$, denoted by $v\in\hat{\partial}f(\bar{x})$, if for every $x\in\dom f$,
			\begin{equation}\label{Formula:frechet subgradient inequality}
				f(x)\geq f(\bar{x})+\ip{v}{x-\bar{x}}+o(\norm{x-\bar{x}}).
			\end{equation}
			\item[(ii)] $v\in\Rn$ is a \textit{limiting subgradient} of $f$ at $\bar{x}\in\dom f$, denoted by $v\in\partial f(\bar{x})$, if
			\begin{equation}\label{Formula:limiting subgraident definition}
				v\in\{v\in\Rn:\exists x_k\xrightarrow[]{f}\bar{x},\exists v_k\in\hat{\partial}f(x_k),v_k\rightarrow v\},
			\end{equation}
			where $x_k\xrightarrow[]{f}\bar{x}\Leftrightarrow x_k\rightarrow\bar{x}\text{ and }f(x_k)\rightarrow f(\bar{x})$. Moreover, we set $\dom\partial f=\{x\in\Rn:\partial f(x)\neq\emptyset\}$. We say that $\bar{x}\in\dom\partial f$ is a stationary point, if $0\in\partial f(\bar{x})$.
		\end{itemize}
	\end{defn}

Recall that the proximal mapping of $f:\Rn\to\overline{\R}$ with parameter $\lambda>0$ is
\[\prox_{\lambda f}(x)=\argmin_{y\in\Rn}\Big\{f(y)+\frac{1}{2\lambda}\norm{x-y}^2\Big\},~\forall x\in\Rn,\]
which is the resolvent $J_{\lambda A}$ with $A=\partial f$ when $f$ is convex; see, e.g.,~\cite{BC}.

For $\eta\in(0,\infty]$, denote by $\Phi_\eta$ the class of functions $\varphi:[0,\eta)\rightarrow\R_+$ satisfying the following conditions: (i) $\varphi(t)$ is right-continuous at $t=0$ with $\varphi(0)=0$; (ii) $\varphi$ is strictly increasing on $[0,\eta)$. The following concept will be the key in our convergence analysis.

\begin{defn}\emph{\cite{wang2020}}\label{Def: g-concave KL} Let $f:\Rn\rightarrow\overline{\mathbb{R}}$
be proper and lsc.  Let $\bar{x}\in\dom\partial f$ and $\mu\in\R$, and
let $V\subseteq\dom\partial f$ be a nonempty subset.
	
(i) We say that $f$ has the pointwise generalized concave Kurdyka-\L ojasiewicz~(KL) property at $\bar{x}\in\dom\partial f$, if there exist neighborhood $U\ni\bar{x}$, $\eta\in(0,\infty]$ and concave $\varphi\in\Phi_\eta$, such that for all $x\in U\cap[0<f-f(\bar{x})<\eta]$,
	\begin{equation}\label{g-concave KL inequality}
		\varphi^\prime_-\big(f(x)-f(\bar{x})\big)\cdot\dist\big(0,\partial f(x)\big)\geq1,
	\end{equation}
where $\varphi_-^\prime$ denotes the left derivative. Moreover, $f$ is a generalized concave KL function if it has the generalized concave KL property at every $x\in\dom\partial f$.
	
(ii) Suppose that $f(x)=\mu$ on $V$. We say $f$ has the setwise\footnote{In the remainder, we shall omit adjectives ``pointwise" and ``setwise" whenever there is no ambiguity.} generalized concave Kurdyka-\L ojasiewicz property on $V$, if there exist $U\supset V$, $\eta\in(0,\infty]$ and concave $\varphi\in\Phi_\eta$ such that for every $x\in U\cap[0<f-\mu<\eta]$,
	\begin{equation}\label{uniform g-concave KL inequality}
		\varphi^\prime_-\big(f(x)-\mu\big)\cdot\dist\big(0,\partial f(x)\big)\geq1.
	\end{equation}
\end{defn}
\begin{remark} The generalized concave KL property reduces to the celebrated concave KL property when concave desingularizing functions $\varphi\in\Phi_\eta$ are continuously differentiable, which has been employed
by a series of papers to establish convergence of proximal-type algorithms in the nonconvex setting; see,~e.g.,~\cite{bolte2018first,Attouch2010,boct2016inertial,attouch2013convergence,
li2016douglas,li2017peaceman,Banert2019,Bolte2014} and the references therein. The generalized concave KL property allows us to describe the optimal~(minimal) concave desingularizing function and
obtain sharp convergence results; see~\cite{wang2020}.
\end{remark}

Many classes of functions satisfy the generalized concave KL property, for instance semialgebraic functions; see, e.g.,~\cite[Corollary 16]{bolte2007} and~\cite[Section 4.3]{Attouch2010}. For more characterizations of the concave KL property, we refer readers to~\cite{bolte2010survey}; see also the fundamental work of \L ojasiewicz~\cite{Lojas1963} and Kurdyka~\cite{Kur98}.

\begin{defn} (i) A set $E\subseteq\mathbb{R}^n$ is called semialgebraic if there exist finitely many polynomials $g_{ij}, h_{ij}:\mathbb{R}^n\rightarrow\mathbb{R}$ such that	
	%\begin{equation*}		
$E=\bigcup_{j=1}^p\bigcap_{i=1}^q\{x\in\mathbb{			R}^n:g_{ij}(x)=0\text{ and }h_{ij}(x)<0\}.	$
	%\end{equation*}	
	
	(ii) A function $f:\mathbb{R}^n\rightarrow\overline{\R}$ is called semialgebraic if its graph $\gph f=\{(x,y)\in\mathbb{R}^{n+1}:f(x)=y\}$ is semialgebraic.
\end{defn}
%\begin{example} The following sets are semialgebraic:
%	Let $A\in\R^{m\times n}$ and $b\in\R^m$. Define $C=\{x\in\Rn:Ax=b\}
%\end{example}
\begin{fact}\label{Fact:semi-algebraic functions are KL} \emph{\cite[Corollary 16]{bolte2007}} Let $f:\mathbb{R}^n\rightarrow\overline{\R}$ be a proper and lsc function and let $\bar{x}\in\dom\partial f$. If $f$ is semialgebraic, then it has the concave KL property at $\bar{x}$ with $\varphi(t)=c\cdot t^{1-\theta}$ for some $c>0$ and $\theta\in(0,1)$.
\end{fact}
%\begin{defn} (i) A set $E\subseteq\mathbb{R}^n$ is called \textit{semialgebraic} if there exist finitely many polynomials $g_{ij}, h_{ij}:\mathbb{R}^n\rightarrow\mathbb{R}$ such that	
%	\begin{equation*}		E=\bigcup_{j=1}^p\bigcap_{i=1}^q\{x\in\mathbb{			R}^n:g_{ij}(x)=0\text{ and }h_{ij}(x)<0\}.	
%	\end{equation*}		(ii) A function $f:\mathbb{R}^n\rightarrow\overline{\R}$ is called \textit{semialgebraic} if its graph		\begin{equation*}		\gph f=\{(x,y)\in\mathbb{R}^{n+1}:f(x)=y\}	\end{equation*}	is semialgebraic.\end{defn}
%\begin{fact}\emph{\cite[Corollary 16]{bolte2007}}\label{Fact:semi-algebraic functions are KL} Let $f:\mathbb{R}^n\rightarrow\overline{\R}$ be a proper and lsc function and let $\bar{x}\in\dom\partial f$. If $f$ is semialgebraic, then it has the concave KL property at $\bar{x}$ with $\varphi(t)=c\cdot t^{1-\theta}$ for some $c>0$ and $\theta\in(0,1)$.\end{fact}

\section{Forward-reflected-backward splitting method for nonconvex problems}\label{section: Forward-reflected-backward splitting method for nonconvex problem}
In this section, we prove convergence results of the FRB splitting method in the nonconvex setting. For
convex minimization problems, the iteration scheme~(\ref{FRB iteration operator form}) specifies to
\begin{align*}
x_{k+1}&=\prox_{\lambda f}\big(y_k-\lambda\nabla g(x_k) \big),
%&=\argmin_{x}\Big\{f(x)+\frac{1}{2\lambda}\norm{x-x_k+2\lambda\nabla g(x_k)-\lambda\nabla g(x_{k-1})}^2\Big\}\\
%&=\argmin_x\Big\{f(x)+\ip{x-y_k}{\nabla g(x_k)}+\frac{1}{2\lambda}\norm{x-y_k}^2\Big\},
\end{align*}
where $y_k=x_k+\lambda\big(\nabla g(x_{k-1})-\nabla g(x_k)\big)$ is an intermediate variable introduced for simplicity.
Since proximal mapping $\prox_{\lambda f}$ might not
be single-valued without convexity of
$f$, we formulate the FRB scheme for solving~(\ref{optimization prblem}) as follows:
\begin{center}
	\fbox{\parbox{\textwidth}
	{\textbf{Forward-reflected-backward (FRB) splitting method}\\
	1. Initialization: Pick $x_{-1}, x_0\in\Rn$ and real number $\lambda>0$.\\
	2. For $k\in\N$, compute
    \begin{align}
	&y_k=x_k+\lambda\big(\nabla g(x_{k-1})-\nabla g(x_k) \big),\label{Algorithm: y_k}\\
	&x_{k+1}\in\prox_{\lambda f}\big(y_k-\lambda\nabla g(x_k) \big).
	\end{align}}}
\end{center}

Recall that a proper, lsc function $f:\Rn\to\overline{\R}$ is prox-bounded if there exists $\lambda>0$ such that $f+\frac{1}{2\lambda}\norm{\cdot}^2$ bounded below; see, e.g.,~\cite[Exercise 1.24]{rockwets}. The maximum of such $\lambda$ is called the proximal threshold, denoted by $\lambda_f$. In the remainder of this paper, we shall
assume that
\begin{itemize}
	\item $f:\Rn\to\overline{\R}$ is proper lsc and prox-bounded with threshold $\lambda_f>0$.
	\item $g:\Rn\to\R$ has a Lipschitz continuous gradient with constant $L>0$.
\end{itemize}

%Fixed $u,v\in\Rn$.
For $0<\lambda<\lambda_f$, \cite[Theorem 1.25]{rockwets} shows that $(\forall x\in\Rn)\ \prox_{\lambda f}(x)$ is nonempty compact,
%the function $x\mapsto f(x)+\ip{x-u}{v}+\frac{1}{2\lambda}\norm{x-u}^2$ is clearly proper, lsc and coercive,  in which case
so the above scheme is well-defined. Moreover, an easy calculation yields that
\begin{equation}\label{Algorithm: x_k}
	\prox_{\lambda f}\big(y_k-\lambda\nabla g(x_k) \big)=\argmin_{x\in\Rn}\Big\{f(x)+\ip{x-y_k}{\nabla g(x_k)}+\frac{1}{2\lambda}\norm{x-y_k}^2\Big\}.
\end{equation}

\subsection{Basic properties}
After formulating the iterative scheme, we now investigate its basic properties. We begin with
 a merit function for the FRB splitting method, which will play a central role in our analysis.
 Such functions of various forms appear frequently in the convergence analysis of splitting methods in nonconvex settings; see, e.g.,~\cite{li2016douglas,li2017peaceman,boct2016inertial}.

\begin{defn}[FRB merit function]
 Let $\lambda>0$. Define the FRB merit function $H:\Rn\times\Rn\to\overline{\R}$ by
	\begin{equation}\label{def: FRB merit function}
		H(x,y)=f(x)+g(x)+\frac{1}{4\lambda}\norm{x-y}^2.
	\end{equation}
\end{defn}
We use the Euclidean norm for $\Rn\times\Rn$, i.e., $\norm{(x,y)}=\sqrt{\norm{x}^2+\norm{y}^2}$.

The FRB merit function has decreasing property under suitable step-size assumption.

\begin{theorem}\label{thm: merit function descent} Let $(x_k)_{k\in\N^*}$ be a sequence generated by FRB and define $z_k=(x_{k+1},x_k)$ for $k\in\N^*$. Assume that $0<\lambda<\min\big\{\frac{1}{4L},\lambda_f\big\}$ and $\inf(f+g)>-\infty$. Let $M_1=\frac{1}{4\lambda}-L$. Then the following statements hold:

(i)	For $k\in\N$, we have
	\begin{align}\label{formula: FRB merit decreases}
		M_1\norm{z_k-z_{k-1}}^2\leq H(z_{k-1})-H(z_k),
	\end{align}
which means that $H(z_k)\leq H(z_{k-1})$. Hence, the sequence $\big(H(z_k)\big)_{k\in\N^*}$ is convergent.
	
(ii) $\displaystyle\sum_{k=0}^\infty\norm{z_k-z_{k-1}}^2<\infty$. Then $\displaystyle\lim_{k\to\infty}\norm{z_k-z_{k-1}}=0$.

\end{theorem}
\proof (i) The FRB iterative scheme and~(\ref{Algorithm: x_k}) imply that
\begin{equation}\label{dd}
	f(x_{k+1})\leq f(x_{k})+\ip{\nabla g(x_k)}{x_{k}-x_{k+1}}+\frac{1}{2\lambda}\norm{x_{k}-y_k}^2-\frac{1}{2\lambda}\norm{x_{k+1}-y_k}^2.
\end{equation}
Invoking the descent lemma~\cite[Theorem 18.15]{BC} to $g$, we have
\begin{equation}\label{ss}
	g(x_{k+1})\leq g(x_k)+\ip{\nabla g(x_k)}{x_{k+1}-x_k}+\frac{L}{2}\norm{x_{k+1}-x_k}^2.
\end{equation}
Summing (\ref{dd}) and (\ref{ss}) yields
\begin{align*}
&~~~~f(x_{k+1})+g(x_{k+1})\leq f(x_k)+g(x_k)+\frac{L}{2}\norm{x_{k+1}-x_k}^2+\frac{1}{2\lambda}\norm{x_{k}-y_k}^2-\frac{1}{2\lambda}\norm{x_{k+1}-y_k}^2\\
	&=f(x_k)+g(x_k)+\left(\frac{L}{2}-\frac{1}{2\lambda}\right)\norm{x_{k+1}-x_k}^2+\frac{1}{\lambda}\ip{x_k-x_{k+1}}{x_k-y_k}\\
	&\leq f(x_k)+g(x_k)+\left(\frac{L}{2}-\frac{1}{2\lambda}\right)\norm{x_{k+1}-x_k}^2+\norm{x_k-x_{k+1}}\norm{\nabla g(x_{k-1})-\nabla g(x_k)}\\
	&\leq f(x_k)+g(x_k)+\left(L-\frac{1}{2\lambda}\right)\norm{x_{k+1}-x_k}^2+\frac{L}{2}\norm{x_k-x_{k-1}}^2,
\end{align*}
where the first equality follows from the identity $\norm{x_k-y_k}^2=\norm{x_{k+1}-y_k}^2-\norm{x_{k+1}-x_k}^2-2\ip{x_{k+1}-x_k}{x_k-y_k}$, and the second inequality is implied by (\ref{Algorithm: y_k}) and the Cauchy-Schwarz inequality. Moreover, the third inequality holds by using the inequality $a^2+b^2\geq 2ab$ and the Lipschitz continuity of~$\nabla g$. Recall the product norm identity $\norm{z_k-z_{k-1}}^2=\norm{x_{k+1}-x_k}^2+\norm{x_k-x_{k-1}}^2$. It then follows from the inequality above that
\begin{align*}
f(x_{k+1})+g(x_{k+1})&\leq f(x_k)+g(x_k)+\left(L-\frac{1}{2\lambda}\right)(\norm{z_k-z_{k-1}}^2-\norm{x_k-x_{k-1}}^2)+\frac{L}{2}\norm{x_k-x_{k-1}}^2\\
&\leq f(x_k)+g(x_k)+\left(L-\frac{1}{2\lambda}\right)\norm{z_k-z_{k-1}}^2+\frac{1}{2\lambda}\norm{x_k-x_{k-1}}^2.
\end{align*}
Now, adding the term $\frac{1}{4\lambda}\norm{x_{k+1}-x_{k}}^2$ to both sides of the above inequality and applying again the product norm identity, one gets
\begin{align*}
H(z_k)&=f(x_{k+1})+g(x_{k+1})+\frac{1}{4\lambda}\norm{x_{k+1}-x_k}^2\nonumber \\
&\leq f(x_k)+g(x_k)+\left(L-\frac{1}{2\lambda}\right)\norm{z_k-z_{k-1}}^2+\frac{1}{2\lambda}\norm{x_k-x_{k-1}}^2\nonumber \\
&+\frac{1}{4\lambda}\norm{z_k-z_{k-1}}^2-\frac{1}{4\lambda}\norm{x_k-x_{k-1}}^2\nonumber \\
&=f(x_k)+g(x_k)+\frac{1}{4\lambda}\norm{x_k-x_{k-1}}^2+\left(L-\frac{1}{4\lambda}\right)\norm{z_k-z_{k-1}}^2\\
&=H(z_{k-1})+\left(L-\frac{1}{4\lambda}\right)\norm{z_k-z_{k-1}}^2,
\end{align*}
which proves (\ref{formula: FRB merit decreases}). By assumption $L-\frac{1}{4\lambda}<0$, which means that $H(z_k)\leq H(z_{k-1})$. The convergence result then follows from the assumption that $\inf (f+g)>-\infty$.

(ii) For integer $l>0$, summing (\ref{formula: FRB merit decreases}) from $k=0$ to $k=l-1$ yields
\[ \sum_{k=0}^{l-1}\norm{z_k-z_{k-1}}^2\leq M_1^{-1}\big(H(z_{-1})-H(z_l)\big)\leq M_1^{-1}\big(H(z_{-1})-\inf(f+g)\big).\]
Passing to the limit, one concludes that $\displaystyle\sum_{k=0}^\infty\norm{z_k-z_{k-1}}^2<\infty$. \endproof
\begin{remark} (Comparison to a known result) Bo{\c{t}} and Csetnek~\cite{boct2016inertial} also utilized $H(x,y)$ as the merit function for Tseng's method, and proved its decreasing property~\cite[Lemma 3.2]{boct2016inertial}. To see the difference between their result and ours, note that the inequality~(\ref{formula: FRB merit decreases}) asserts that $H$ decreases with respect to a single sequence~$(x_k)$ generated by FRB, while~\cite[Lemma 3.2]{boct2016inertial} concerns the decreasing property of $H$ with respect to two difference sequences produced by Tseng's method.
\end{remark}

Next we estimate the lower bound on the gap between iterates. To this end, a lemma helps. % through the limiting subgradients of the FRB merit function.
\begin{lemma}\label{lem: merit function subdifferential} For every $(x,y)\in\Rn\times\Rn$, we have
$$\partial H(x,y)=\left\{\partial f(x)+\nabla g(x)+\frac{1}{2\lambda}(x-y)\right\}\times
\left\{\frac{1}{2\lambda}(y-x)\right\}.$$
\end{lemma}
\proof Apply the subdifferential sum and separable sum rules~\cite[Excercise 8.8, Proposition 10.5]{rockwets}.\endproof

\begin{theorem}\label{thm: subgradient lower bound} Let $(x_k)_{k\in\N^*}$ be a sequence generated by FRB and define $z_k=(x_{k+1},x_k)$ for $k\in\N^*$. Let $M_2=\sqrt{2}\left(L+\frac{2}{\lambda}\right)$ and define for $k\in\N$
\begin{align*}
&p_{k+1}=\frac{1}{\lambda}(x_k-x_{k+1})+\nabla g(x_{k+1})-2\nabla g(x_k)+\nabla g(x_{k-1}),\\
&A_{k+1}=p_{k+1}+\frac{1}{2\lambda}(x_{k+1}-x_k),~B_{k+1}=\frac{1}{2\lambda}(x_k-x_{k+1}).
\end{align*}
Then $\left(A_{k+1},B_{k+1}\right)\in\partial H(z_k)$ and $\norm{\left(A_{k+1},B_{k+1}\right)}\leq M_2\norm{z_k-z_{k-1}}$.
% $k\in\N$.
\end{theorem}
\proof The FRB scheme and~(\ref{Algorithm: x_k}) imply that $0\in\partial f(x_{k+1})+\nabla g(x_k)+\frac{1}{\lambda}(x_{k+1}-y_k)$. Equivalently,
$\frac{1}{\lambda}(x_k-x_{k-1})+\nabla g(x_{k-1})-2\nabla g(x_k)\in\partial f(x_{k+1})$, which together with the subdifferential sum rule~\cite[Excercise 8.8]{rockwets} implies that
\begin{align*}
p_{k+1}\in\partial f(x_{k+1})+\nabla g(x_{k+1}).
%=\partial (f+g)(x_{k+1}).
\end{align*}
Applying Lemma~\ref{lem: merit function subdifferential} with $x=x_{k+1}$ and $y=x_k$ yields \[\partial H(z_k)=\Big\{\partial f(x_{k+1})+\nabla g(x_{k+1})+\frac{1}{2\lambda}(x_{k+1}-x_k)\Big\}
\times\Big\{\frac{1}{2\lambda}(x_k-x_{k+1})\Big\}.\]
It then follows that $\left(A_{k+1},B_{k+1}\right)=\left(p_{k+1}+\frac{1}{2\lambda}(x_{k+1}-x_k),\frac{1}{2\lambda}(x_k-x_{k+1})\right)\in\partial H(z_k)$ and
\begin{align*}
&~~~\norm{\left(A_{k+1},B_{k+1}\right)}\leq\norm{A_{k+1}}+\norm{B_{k+1}}=\norm{p_{k+1}+\frac{1}{2\lambda}(x_{k+1}-x_k)}+\norm{\frac{1}{2\lambda}(x_{k}-x_{k+1})}\\
&\leq\frac{2}{\lambda}\norm{x_{k+1}-x_k}+\norm{\nabla g(x_{k+1})-\nabla g(x_k)}+\norm{\nabla g(x_k)-\nabla g(x_{k-1})}\\
&\leq\left(\frac{2}{\lambda}+L\right)\norm{x_{k+1}-x_k}+\left(\frac{2}{\lambda}+L\right)\norm{x_{k}-x_{k-1}}\leq M_2\norm{z_k-z_{k-1}},
\end{align*}
where the third inequality follows from the Lipschitz continuity of $\nabla g$ with modulus $L$. \endproof

We now connect the properties of FRB merit function
with the actual objective of~(\ref{optimization prblem}).
Denote by $\omega(z_{-1})$ the set of limit points of $(z_k)_{k\in\N^*}$.

\begin{theorem}\label{thm: function values convergence} Let $(x_k)_{k\in\N^*}$ be a sequence generated by FRB and define $z_k=(x_{k+1},x_k)$ for $k\in\N^*$. Assume that conditions in Theorem~\ref{thm: merit function descent} are satisfied and $(z_k)_{k\in\N^*}$ is bounded. Suppose that a subsequence $(z_{k_l})_{l\in\N}$ of $(z_k)_{k\in\N^*}$ converges to some $z^*=(x^*,y^*)$ as $l\to\infty$. Then the following statements hold:

(i) $\displaystyle\lim_{l\to\infty} H(z_{k_l})=f(x^*)+g(x^*)=F(x^*)$. In fact, $\displaystyle\lim_{k\to\infty} H(z_{k})=F(x^*)$.

(ii) We have $x^*=y^*$ and $0\in\partial H(x^*,y^*)$, which implies $0\in
\partial F(x^*)=\partial f(x^*)+\nabla g(x^*)$.

(iii) The set $\omega(z_{-1})$ is nonempty, compact and connected, on which the FRB merit function $H$ is finite and constant. Moreover, we have~$\displaystyle\lim_{k\to\infty}\dist(z_k,\omega(z_{-1}))=0$.
\end{theorem}
\proof (i) By the identity~(\ref{Algorithm: x_k}), for every $k\in\N^*$
\begin{align*}
f(x_{k+1})+\ip{x_{k+1}-y_k}{\nabla g(x_k)}+\frac{1}{2\lambda}\norm{x_{k+1}-y_k}^2\leq f(x^*)+\ip{x^*-y_k}{\nabla g(x_k)}+\frac{1}{2\lambda}\norm{x^*-y_k}^2.
\end{align*}
Combining the identity $\norm{x^*-y_k}^2=\norm{x^*-x_{k+1}}^2+\norm{x_{k+1}-y_k}^2+2\ip{x^*-x_{k+1}}{x_{k+1}-y_k}$ with the above inequality and using the definition $y_k=x_k+\lambda\big(\nabla g(x_{k-1})-\nabla g(x_k) \big)$, one gets
\begin{align*}
f(x_{k+1})&\leq f(x^*)+\ip{x^*-x_{k+1}}{\nabla g(x_k)}+\frac{1}{2\lambda}\norm{x^*-x_{k+1}}^2+\frac{1}{\lambda}\ip{x^*-x_{k+1}}{x_{k+1}-y_k}\\
&\leq f(x^*)+\ip{x^*-x_{k+1}}{2\nabla g(x_k)-\nabla g(x_{k-1})}+\frac{1}{2\lambda}\norm{x^*-x_{k+1}}^2+\frac{1}{\lambda}\ip{x^*-x_{k+1}}{x_{k+1}-x_k}.
%&\leq f(x^*)+\norm{x^*-x_{k+1}}\norm{2\nabla g(x_k)-\nabla g(x_{k-1})}+\frac{1}{2\lambda}\norm{x^*-x_{k+1}}^2+\frac{1}{\lambda}\norm{x^*-x_{k+1}}\norm{x_{k+1}-x_k}
\end{align*}
Replacing $k$ by $k_l$, passing to the limit and applying Theorem~\ref{thm: merit function descent}(ii), we have
\[\limsup_{l\to\infty}f(x_{k_l+1})\leq f(x^*)\leq\liminf_{l\to\infty}f(x_{k_l+1}),\]
where the last inequality is implied by the lower semicontinuity of $f$, which implies that $\displaystyle\lim_{l\to\infty}f(x_{k_l+1})=f(x^*)$. Hence we conclude that $$\displaystyle\lim_{l\to\infty}H(z_{k_l})=\lim_{l\to\infty}\left(f(x_{k_l+1})+g(x_{k_l+1})+\frac{1}{4\lambda}\norm{x_{k_l+1}-x_{k_l}}^2\right)=f(x^*)+g(x^*),$$
where the last equality holds because of Theorem~\ref{thm: merit function descent}(ii). Moreover, Theorem~\ref{thm: merit function descent}(i) states that the function value sequence $\big(H(z_k)\big)_{k\in\N^*}$ converges, from which the rest of the statement readily follows.

(ii) By assumption, we have $x_{k_l+1}\to x^*$ and $x_{k_l}\to y^*$ as $l\to\infty$. Therefore $0\leq\norm{y^*-x^*}\leq\norm{y^*-x_{k_l}}+\norm{x_{k_l}-x_{k_l+1}}+\norm{x_{k_l+1}-x^*}\to0$ as $l\to\infty$, which means that $x^*=y^*$. Combining Theorem~\ref{thm: merit function descent}, Theorem~\ref{thm: subgradient lower bound}, statement(i) and the outer semicontinuity of $z\mapsto\partial H(z)$, one gets \[0\in\partial H(x^*,y^*)=\partial H(x^*,x^*)=\big\{\partial f(x^*)+\nabla g(x^*)\big\}\times\big\{0\big\}\Rightarrow 0\in\partial f(x^*)+ \nabla g(x^*),\]
as desired.

(iii) Given that $(z_k)_{k\in\N^*}$ is bounded, it is easy to see that $\omega(z_{-1})=\cap_{l\in\N}\overline{\cup_{k\geq l}z_k}$ is a nonempty compact set. By using Theorem~\ref{thm: merit function descent} and a result by Ostrowski~\cite[Theorem 1.49]{BC}, the set $\omega(z_{-1})$ is connected. Pick $\tilde{z}\in\omega(z_{-1})$ and suppose that $z_{k_q}\to\tilde{z}$ as $q\to\infty$. Then statement(i) implies that $H(\tilde{z})=\displaystyle\lim_{q\to\infty}H(z_{k_q})=H(z^*)$. Finally, by the definition of $\omega(z_{-1})$, we have $\displaystyle\lim_{k\to\infty}\dist(z_k,\omega(z_{-1}))=0$.\endproof

Theorem~\ref{thm: function values convergence} requires the sequence $(z_k)_{k\in\N^*}$ to be bounded. The result below provides a sufficient condition to such assumption.

\begin{theorem}\label{thm: bounded sequence} Let $(x_k)_{k\in\N^*}$ be a sequence generated by FRB and define $z_k=(x_{k+1},x_k)$ for $k\in\N^*$. Assume that conditions in Theorem~\ref{thm: merit function descent} are satisfied. If $f+g$ is coercive (or level bounded), then the
sequence $(x_k)_{k\in\N^*}$ is bounded, so is $(z_k)_{k\in\N^*}$.
\end{theorem}
\proof Theorem~\ref{thm: merit function descent}(i) implies that we have
\begin{align*}
	H(x_1,x_0)&\geq f(x_{k+1})+g(x_{k+1})+\frac{1}{4\lambda}\norm{x_{k+1}-x_k}^2\geq f(x_{k+1})+g(x_{k+1}).
\end{align*}
Suppose that $(x_k)_{k\in\N^*}$ was unbounded. Then we would have a contradiction by the coercivity or
level boundedness of $f+g$.
\endproof
\subsection{Convergence under the generalized concave KL property}

Following basic properties of the FRB method, we now present the main convergence result under the generalized concave KL property. For $\varepsilon>0$ and nonempty set $\Omega\subseteq\Rn$, we define $\Omega_\varepsilon=\{x\in\Rn:\dist(x,\Omega)<\varepsilon\}$. The following lemma will be useful soon.

\begin{lemma}\label{lemma: Uniformize the g-concave KL}\emph{\cite[Lemma 4.4]{wang2020}} Let $f:\Rn\rightarrow\overline{\R}$ be proper lsc and let $\mu\in\R$. Let $\Omega\subseteq\dom\partial f$ be a nonempty compact set on which $f(x)=\mu$ for all $x\in\Omega$. The following statements hold:
	
(i) Suppose that $f$ satisfies the pointwise generalized concave KL property at each $x\in\Omega$. Then there exist $\varepsilon>0,\eta\in(0,\infty]$ and $\varphi(t)\in\Phi_\eta$ such that $f$ has the setwise generalized concave KL property on $\Omega$ with respect to $U=\Omega_\varepsilon$, $\eta$ and $\varphi$.

(ii) Set $U=\Omega_\varepsilon$ and define $h:(0,\eta)\rightarrow\R_+$ by \[h(s)=\sup\big\{\dist^{-1}\big(0,\partial f(x)\big):x\in U\cap[0<f-\mu<\eta],s\leq f(x)-\mu\big\}.\]
Then the function $\tilde{\varphi}:[0,\eta)\rightarrow\R_+$,
\[t\mapsto\int_0^th(s)ds,~\forall t\in(0,\eta),\]
and $\tilde{\varphi}(0)=0$, is well-defined and belongs to $\Phi_\eta$. The function $f$ has the setwise generalized concave KL property on $\Omega$ with respect to $U$, $\eta$ and $\tilde{\varphi}$. Moreover,
\[\tilde{\varphi}=\inf\big\{\varphi\in\Phi_\eta:\text{$\varphi$ is a concave desingularizing function of $f$ on $\Omega$ with respect to $U$ and $\eta$}\big\}.\]
We say $\tilde{\varphi}$ is the exact modulus of the setwise generalized concave KL property of $f$ on $\Omega$ with respect to $U$ and $\eta$.
\end{lemma}

We are now ready for the main results. Our methodology is akin to the concave KL convergence mechanism employed by a vast amount of literature; see, e.g.,~\cite{bolte2018first,Attouch2010,boct2016inertial,attouch2013convergence,li2016douglas,li2017peaceman,Banert2019}, but makes use of the generalized concave KL property and the associated exact modulus, which guarantees the sharpness of our results; see Remark~\ref{rem: sharp result}.

\begin{theorem}[Global convergence of FRB]\label{thm: finite length}
 Let $(x_k)_{k\in\N^*}$ be a sequence generated by FRB and define $z_k=(x_{k+1},x_k)$ for $k\in\N^*$. Assume that $(z_k)_{k\in\N^*}$ is bounded, $\inf(f+g)>-\infty$, and $0<\lambda<\min\big\{\frac{1}{4L},\lambda_f\big \}$. Suppose that the FRB merit function $H(x,y)$ has the generalized concave KL property on $\omega(z_{-1})$. Then the following statements hold:
	
(i) The sequence $(z_k)_{k\in\N^*}$ is Cauchy and has finite length. To be specific, there exist $M>0$, $k_0\in\N$, $\varepsilon>0$ and $\eta\in(0,\infty]$ such that for $i\geq k_0+1$
\begin{equation}\label{formula: finite length}
	\sum_{k=i}^\infty\norm{z_{k+1}-z_k}\leq\norm{z_{i}-z_{i-1}}+M\tilde{\varphi}\left( H(z_{i})-H(z^*)\right).
\end{equation}
where $\tilde{\varphi}\in\Phi_\eta$ is the exact modulus associated with the setwise generalized concave KL property of $H$ on $\omega(z_{-1})$ with respect to $\varepsilon$ and $\eta$.

(ii) The sequence $(x_k)_{k\in\N^*}$ has finite length and converges to some $x^*$ with $0\in\partial F(x^*)$.
\end{theorem}
\proof (i) By the boundedness assumption, assume without loss of generality that $z_k\to z^*=(x^*,y^*)$ for some $z^*\in\Rn\times\Rn$. Then Theorem~\ref{thm: function values convergence} implies that $x^*=y^*$ and $H(z_k)\to F(x^*)=H(z^*)$ as $k\to\infty$. Recall from Theorem~\ref{thm: merit function descent} that we have $H(z_{k+1})\leq H(z_k)$ for $k\in\N^*$, therefore one needs to consider two cases.

Case 1: Suppose that there exists $k_0$ such that $H(z^*)=H(z_{k_0})$. Then Theorem~\ref{thm: merit function descent}(i) implies that $z_{k_0+1}=z_{k_0}$ and $x_{k_0+1}=x_{k_0}$. The desired results then follows from a simple induction.

Case 2: Assume that $H(z^*)<H(z_k)$ for every $k$. By assumption and Lemma~\ref{lemma: Uniformize the g-concave KL}, there exist $\varepsilon>0$ and $\eta\in(0,\infty]$ such that the FRB merit function $H$ has the setwise generalized concave KL property on $\omega(z_{-1})$ with respect to $\varepsilon>0$ and $\eta>0$ and the associated exact modulus $\tilde{\varphi}$. On one hand, by the fact that $H(z_k)\to H(z^*)$, there exists $k_1$ such that $z_k\in[0<H-H(z^*)<\eta]$ for $k>k_1$. On the other hand, Theorem~\ref{thm: function values convergence}(iii) implies that there exists $k_2$ such that $\dist(z_k,\omega(z_{-1}))<\varepsilon$ for all $k>k_2$. Put $k_0=\max(k_1,k_2)$. Then for $k>k_0$, we have $z_k\in \omega(z_{-1})_\varepsilon\cap[0<H-H(z^*)<\eta]$. %By using Theorem~\ref{thm: function values convergence}(iii) and Lemma~\ref{lemma: Uniformize the g-concave KL}, we can see that the FRB merit function $H$ has the setwise generalized concave KL property on $\omega(z_{-1})$ with respect to some $\varepsilon>0$ and $\eta>0$ and the associated exact modulus $\tilde{\varphi}$.
Hence for $k>k_0$
\begin{equation*}
	(\tilde{\varphi})_-^\prime\left(H(z_k)-H(z^*)\right)\cdot\dist(0,\partial H(z_k))\geq1.
\end{equation*}
Invoking Theorem~\ref{thm: subgradient lower bound} yields
\begin{equation}\label{d}
M_2(\tilde{\varphi})_-^\prime\left(H(z_k)-H(z^*)\right)\norm{z_k-z_{k-1}}\geq1.
\end{equation}

For simplicity, define for $k>l$ $$\Delta_{k,k+1}=\tilde{\varphi}\left(H(z_k)-H(z^*)\right)-\tilde{\varphi}\left(H(z_{k+1})-H(z^*)\right).$$ By the concavity of $\tilde{\varphi}$ and~(\ref{d}), one has
\begin{align}\label{cc}
\Delta_{k,k+1}&=\tilde{\varphi}\left(H(z_k)-H(z^*)\right)-\tilde{\varphi}\left(H(z_{k+1})-H(z^*)\right)\nonumber \\
&\geq(\tilde{\varphi})_-^\prime\left(H(z_k)-H(z^*)\right)\cdot[H(z_k)-H(z_{k+1})]\nonumber\\
&\geq\frac{H(z_k)-H(z_{k+1})}{M_2\norm{z_k-z_{k-1}}}.
\end{align}
Applying Theorem~\ref{thm: merit function descent} to~(\ref{cc}) implies that
\begin{equation*}
	\Delta_{k,k+1}\geq\frac{\norm{z_{k+1}-z_k}^2}{\norm{z_k-z_{k-1}}}\cdot\frac{M_1}{M_2}.
\end{equation*}
Put $M=\frac{M_2}{M_1}$. Note that $a^2+b^2\geq2ab$ for $a,b\geq0$. Then the above inequality gives
\begin{equation}\label{xx}
	2\norm{z_{k+1}-z_k}\leq 2\sqrt{M\Delta_{k,k+1}\norm{z_k-z_{k-1}}}\leq M\Delta_{k,k+1}+\norm{z_k-z_{k-1}}.
\end{equation}
Pick $i\geq k_0+1$. Summing~(\ref{xx}) from $i$ to an arbitrary $j>i$ yields
\begin{align*}
&~~~~2\sum_{k= i }^ j\norm{z_{k+1}-z_k}\leq\sum_{k= i }^ j\norm{z_k-z_{k-1}}+M\sum_{k= i }^ j\Delta_{k,k+1}\\
&\leq\sum_{k=i}^ j\norm{z_{k+1}-z_k}+\norm{z_{ i }-z_{i-1}}+M\tilde{\varphi}\left( H(z_{ i })-H(z^*)\right)-M\tilde{\varphi}\left( H(z_{ j+1})-H(z^*)\right)\\
&\leq \sum_{k= i }^ j\norm{z_{k+1}-z_k}+\norm{z_{ i }-z_{i-1}}+M\tilde{\varphi}\left( H(z_{ i })-H(z^*)\right),
\end{align*}
where the second inequality is implied by the definition of $\Delta_{k,k+1}$, implying
\[\sum_{k= i }^ j\norm{z_{k+1}-z_k}\leq\norm{z_i-z_{i-1}}+M\tilde{\varphi}\left(H(z_i)-H(z^*)\right),\]
from which~(\ref{formula: finite length}) readily follows. Moreover, one gets from the above inequality that for $i\geq k_0+1$ and $j>i$, \[\norm{z_{i+j}-z_i}\leq\sum_{k=i}^{j-1}\norm{z_{k+1}-z_k}\leq\norm{z_i-z_{i-1}}+M\tilde{\varphi}(H(z_i)-H(z^*)).\]
Recall from Theorem~\ref{thm: function values convergence}(i)\&(ii) that $H(z_i)\to H(z^*)$ and from Theorem~\ref{thm: merit function descent}(ii) that $\norm{z_i-z_{i-1}}\to0$ as $i\to\infty$. Passing to the limit, one concludes that $(z_k)_{k\in\N^*}$ is Cauchy.

(ii) The statement follows from the definition of $(z_k)_{k\in\N^*}$ and Theorem~\ref{thm: function values convergence}(ii).\endproof
\begin{remark}\label{rem: sharp result}(i) Note that iterates distance was shown to be only square-summable in the original FRB paper~\cite[Theorem 2.5]{malitsky2020FRB}. Therefore the finite length property is even new in the convex setting.
	
(ii) Unlike the usual concave KL convergence analysis, our approach uses the generalized concave KL property and the associated exact modulus to describe the sharpest upper bound of $\sum_{k=-1}^\infty\norm{z_{k+1}-z_k}$. To see this, note that the usual analysis would yield a similar bound as~(\ref{formula: finite length}) with $\tilde{\varphi}$ replaced by a concave desingularizing function associated with the concave KL property of $H$. Lemma~\ref{lemma: Uniformize the g-concave KL} states that $\tilde{\varphi}$ is the infimum of all associated concave desingularizing functions. Hence the upper bound~(\ref{formula: finite length}) is the sharpest; see also~\cite{wang2020} for a similar sharp result of the celebrated PALM algorithm.
\end{remark}

A key assumption of Theorem~\ref{thm: finite length} is the concave KL property of the FRB merit function $H$ on $\omega(z_{-1})$.  The class of semialgebraic functions provides rich examples of functions satisfying such an assumption.
%Another celebrated result states that semialgebraic functions have the concave KL property.
%\begin{proposition} Let $f:\Rn\to\overline{\R}$ be proper and lsc. Define $H(x,y)=f(x)+r\norm{x-y}^2$ for $r>0$. If $f$ is semialgebraic, then so is $H$. Consequently, $H$ is a concave KL function.\end{proposition}
%\proof Recall from~\cite[Section 4.3]{Attouch2010} that the class of semialgebraic functions is closed under summation and notice that the quadratic function $(x,y)\mapsto r\norm{x-y}^2$ is semialgebraic. Applying Fact~\ref{Fact:semi-algebraic functions are KL} completes the proof.\endproof

\begin{corollary}\label{cor: semialgebraic convergence} Let $(x_k)_{k\in\N^*}$ be a sequence generated by FRB. Assume that $(x_k)_{k\in\N^*}$ is bounded, $\inf(f+g)>-\infty$, and $0<\lambda<\min\big\{\frac{1}{4L},\lambda_f \big\}$.  Suppose further that $f$ and $g$ are both semialgebraic functions.  Then $(x_k)_{k\in\N^*}$ converges to some $x^*$ with $0\in\partial F(x^*)$ and has finite length property.
\end{corollary}
\proof Recall from~\cite[Section 4.3]{Attouch2010} that the class of semialgebraic functions is closed under summation and notice that the quadratic function $(x,y)\mapsto \frac{1}{4\lambda}\norm{x-y}^2$ is semialgebraic. Then Fact~\ref{Fact:semi-algebraic functions are KL} implies that the FRB merit function $H$ is concave KL. Applying Theorem~\ref{thm: finite length} then completes the proof.
\endproof

Assuming that the FRB merit function admits KL exponent $\theta\in[0,1)$, we establish the following convergence rates result. Our analysis is standard and follows from the usual KL convergence rate methodology; see, e.g.,~\cite[Theorem 2]{attouch2009convergence},~\cite[Remark 6]{Bolte2014} and~\cite[Lemma 4]{Banert2019}.  We provide a proof here for the sake of completeness.

\begin{theorem}[Convergence rate of FRB]\label{thm: rate}
Let $(x_k)_{k\in\N^*}$ be a sequence generated by FRB and define $z_k=(x_{k+1},x_k)$ for $k\in\N^*$. Assume that $(z_k)_{k\in\N^*}$ is bounded, $\inf(f+g)>-\infty$, and $0<\lambda<\min\big\{\frac{1}{4L},\lambda_f\big\}$. Let $M,k_0$ and $x^*$ be those given in Theorem~\ref{thm: finite length}.  Suppose that the FRB merit function $H(x,y)$ has KL exponent $\theta\in[0,1)$ at $(x^*,x^*)$. %Suppose further the associated desingularizing function of $H$ has the form $\varphi(t)=c\cdot t^{1-\theta}$ for some $c>0$ and $\theta\in[0,1)$.
Then the following statements hold:

(i) If $\theta=0$, then $(x_k)_{k\in\N^*}$ converges to $x^*$ in finite steps.

(ii) If $\theta\in\left(0,\frac{1}{2}\right]$, then there exist $Q_1\in(0,1)$ and $c_1>0$ such that \[\norm{x_k-x^*}\leq Q_1^{k}c_1,\]
 for $k$ sufficiently large.

(iii) If $\theta\in\left(\frac{1}{2},1\right)$, then there exists $Q_2>0$ such that
\[\norm{x_k-x^*}\leq Q_2\cdot k^{-\frac{1-\theta}{2\theta-1}},\]
for $k$ sufficiently large.
\end{theorem}
\proof Let $z^*=(x^*,x^*)$. Then Theorem~\ref{thm: finite length} shows that $z_k\to z^*$ as $k\to\infty$. Assume without loss of generality that $H(z^*)=0$ and $H(z_k)>0$ for all $k$. %We shall only consider the case where $H(z_{k+1})<H(z_k)$ for every $k$, otherwise a similar argument as in Theorem~\ref{thm: finite length} shows that $x_k$ converges in finite steps, in which case all desired results are trivial.
Before proving the desired statements, we will first develop several inequalities which will be used later. We have shown in Theorem~\ref{thm: finite length} that for $k>k_0$,
\begin{align*}
1\leq(\tilde{\varphi})_-^\prime\left(H(z_k)\right)\cdot\dist(0,\partial H(z_k))\leq\varphi_-^\prime\left(H(z_k)\right)\cdot\dist(0,\partial H(z_k)),
\end{align*}
which by Theorem~\ref{thm: subgradient lower bound} and our assumption further implies that
\begin{align}\label{x}
1\leq (1-\theta)c\left(H(z_k)\right)^{-\theta}\norm{(A_{k+1},B_{k+1})}\leq(1-\theta) cM_2\left(H(z_k)\right)^{-\theta}\norm{z_k-z_{k-1}}.
\end{align}
Furthermore, define for $k\in\N$
\[\sigma_k=\sum_{i=k}^\infty\norm{z_{i+1}-z_i},\]
which is well-defined due to Theorem~\ref{thm: finite length}. Assume again without loss of generality that $\sigma_{k-1}-\sigma_k=\norm{z_k-z_{k-1}}<1$ for every $k$ (recall Theorem~\ref{thm: merit function descent}). Rearranging~(\ref{x}) gives
\begin{equation}\label{cd}
H(z_k)^{1-\theta}\leq\left[(1-\theta)cM_2\right]^{\frac{1-\theta}{\theta}}(\sigma_{k-1}-\sigma_k)^\frac{1-\theta}{\theta}.
\end{equation}

By~(\ref{formula: finite length}), we have for $k>k_0$
\begin{equation}\label{sigma_k upper bound}
\sigma_k\leq\sigma_{k-1}-\sigma_k+cM(H(z_k))^{1-\theta}\leq\sigma_{k-1}-\sigma_k+C(\sigma_{k-1}-\sigma_k)^{\frac{1-\theta}{\theta}},
\end{equation}
where the last inequality follows from~(\ref{cd}) and $C=cM((1-\theta)cM_2)^{\frac{1-\theta}{\theta}}$. Clearly $\norm{x_k-x^*}\leq\norm{z_k-z^*}\leq\sigma_k$. Hence, in order to prove the desired statements, it suffices to estimate $\sigma_k$.

(i): Let $\theta=0$ and suppose that $(z_k)_{k\in\N^*}$ converges in infinitely many steps. Then~(\ref{x}) means that for every $k>k_0$ %$H(z_k)^\theta\leq(1-\theta)CM_2\norm{z_k-z_{k-1}}$, and $H(z_k)>0$. Hence
$$1\leq cM_2\norm{z_k-z_{k-1}}.$$
Passing to the limit and applying Theorem~\ref{thm: merit function descent}, one gets $1\leq 0$, which is absurd.

(ii) If $\theta\in\left(0,\frac{1}{2}\right]$, then $\frac{1-\theta}{\theta}\geq1$ and $(\sigma_{k-1}-\sigma_k)^{\frac{1-\theta}{\theta}}\leq\sigma_{k-1}-\sigma_k$. Hence~(\ref{sigma_k upper bound}) implies that for $k>k_0$
\begin{equation*}
\sigma_k\leq(C+1)(\sigma_{k-1}-\sigma_k)\Rightarrow\sigma_k\leq\frac{C+1}{C+2}\sigma_{k-1}\leq\left(\frac{C+1}{C+2}\right)^{k-k_0}\sigma_{k_0},
\end{equation*}
from which the desired result readily follows by setting $Q_1=\frac{C+1}{C+2}$ and $c_1=\left(\frac{C+1}{C+2}\right)^{-k_0}\sigma_{k_0}$.

(iii) If $\frac{1}{2}<\theta<1$, then $0<\frac{1-\theta}{\theta}<1$ and $(\sigma_{k-1}-\sigma_k)^{\frac{1-\theta}{\theta}}\geq\sigma_{k-1}-\sigma_k$. It follows from~(\ref{sigma_k upper bound}) that $\sigma_k\leq\left(1+C \right)(\sigma_{k-1}-\sigma_k)^{\frac{1-\theta}{\theta}}$. Define $h(t)=t^{-\frac{\theta}{1-\theta}}$ for $t>0$. Then
\begin{equation}\label{c}
1\leq\left(1+C \right)^{\frac{\theta}{1-\theta}}(\sigma_{k-1}-\sigma_k)h(\sigma_k).
\end{equation}
Let $R>1$ be a fixed real number. Next we consider two cases for $k>k_0$. If $h(\sigma_{k})\leq Rh(\sigma_{k-1})$ then~(\ref{c}) implies that
\begin{align*}
1&\leq R(1+C)^{\frac{\theta}{1-\theta}}(\sigma_{k-1}-\sigma_k)h(\sigma_{k-1})\leq R(1+C)^{\frac{\theta}{1-\theta}}\int_{\sigma_k}^{\sigma_{k-1}}h(t)dt\\
&=R(1+C)^{\frac{\theta}{1-\theta}}\frac{1-\theta}{1-2\theta}\left(\sigma_{k-1}^{\frac{1-2\theta}{1-\theta}}-\sigma_k^{\frac{1-2\theta}{1-\theta}}\right).
\end{align*}
Set $v=\frac{1-2\theta}{1-\theta}<0$. Then the above inequality can be rewritten as
\begin{equation}\label{v1}
\sigma_{k}^v-\sigma_{k-1}^v\geq-\frac{v}{R}(1+C)^{-\frac{\theta}{1-\theta}}.
\end{equation}
If $h(\sigma_k)>Rh(\sigma_{k-1})$ then $\sigma_k<\left(\frac{1}{R}\right)^{\frac{1-\theta}{\theta}}\sigma_{k-1}=q\sigma_{k-1}$, where $q=\left(\frac{1}{R}\right)^{\frac{1-\theta}{\theta} }\in(0,1)$. Hence
\begin{equation}\label{v2}
\sigma_k<q\sigma_{k-1}\Rightarrow\sigma_k^v>q^v\sigma_{k-1}^v\Rightarrow\sigma_k^v-\sigma_{k-1}^v>(q^v-1)\sigma_{k-1}^v>(q^v-1)\sigma_{k_0}^v.
\end{equation}
Set $c_2=\min\big\{(q^v-1)\sigma_{k_0}^v,-\frac{v}{R}(1+C)^{-\frac{\theta}{1-\theta}}\big\}>0$. Combining~(\ref{v1}) and~(\ref{v2}) yields $$\sigma_k^v-\sigma_{k-1}^v\geq c_2,\forall k>k_0.$$
Summing the above inequality from $k_0+1$ to any $k>k_0$, one gets
\begin{equation*}
\sigma_k^v-\sigma_{k_0}^v=\sum_{i=k_0+1}^k(\sigma_i^v-\sigma_{i-1}^v)\geq(k-k_0)c_2\Rightarrow\sigma_k^v>\frac{k}{k_0+1}c_2,
\end{equation*}
where the last inequality holds because $k-k_0\geq\frac{k}{k_0+1}$ for $k>k_0+1$. Hence $$\sigma_k<k^{\frac{1}{v}}\left(\frac{C_2}{k_0+1}\right)^{\frac{1}{v}},$$ from which the desired statement follows by setting $Q_2=\left(\frac{c_2}{k_0+1}\right)^{\frac{1}{v}}$.\endproof

The following corollary asserts that convergence rates can be directly deduced from KL exponent of the objective under certain conditions, which is a consequence of~\cite[Theorem 3.6]{li2018calculus}, Theorems~\ref{thm: function values convergence} and~\ref{thm: rate}.

\begin{corollary}\label{prop: KL exponent of merit function} Let $(x_k)_{k\in\N^*}$ be a sequence generated by FRB and suppose that all conditions of Theorem~\ref{thm: finite length} are satisfied. Let $x^*$ be as in Theorem~\ref{thm: finite length}. Assume further that $F$ has KL exponent $\theta\in[\frac{1}{2},1)$ at $x^*$. Then the following statements hold:
	
	(i) If $\theta=\frac{1}{2}$, then there exist $Q_1\in(0,1)$ and $c_1>0$ such that \[\norm{x_k-x^*}\leq Q_1^{k}c_1,\]
	for $k$ sufficiently large.
	
	(ii) If $\theta\in\left(\frac{1}{2},1\right)$, then there exists $Q_2>0$ such that
	\[\norm{x_k-x^*}\leq Q_2\cdot k^{-\frac{1-\theta}{2\theta-1}},\]
	for $k$ sufficiently large.
\end{corollary}
\proof The FRB merit function $H$ has KL exponent $\theta$ at $(x^*,x^*)$ by our assumption and~\cite[Theorem 3.6]{li2018calculus}. Applying Theorem~\ref{thm: rate} completes the proof.
\endproof
\begin{example}\label{ex: feasibility convergence} Let $C\subseteq\Rn$ be a nonempty, closed and convex set and let $D\subseteq\Rn$ be nonempty and closed. Suppose that $C\cap D\neq\varnothing$ and either $C$ or $D$ is compact. Assume further that both $C$ and $D$ are  semialgebraic.  Consider the minimization problem~(\ref{optimization prblem}) with $f=\delta_D$ and $g=\frac{1}{2}\dist^2(\cdot,C)$. Let $0<\lambda<\frac{1}{4}$, and let $(x_k)_{k\in\N^*}$ be a sequence generated by FRB and $z_k=(x_{k+1},x_k)$ for $k\in\N^*$. Then the following statements hold:

(i) There exists $x^*\in\Rn$ such that $x_k\to x^*$ and $0\in\partial F(x^*)$.

(ii) Suppose in addition that
\begin{equation}\label{CQ}
			N_C(\proj_C(x^*) )\cap\big(- N_D(x^*) \big)=\{0\}.
	\end{equation}
Then $x^*\in C\cap D$. Moreover, there exist $Q_1\in(0,1)$ and $c_1>0$ such that \[\norm{x_k-x^*}\leq Q_1^{k}c_1,\]
for $k$ sufficiently large.
\end{example}
\proof (i) By the compactness assumption, the function $f+g$ is coercive. Hence Theorem~\ref{thm: bounded sequence} implies that $(x_k)_{k\in\N^*}$ is bounded. We assume that sets $C,D$ are semialgebraic, then so are functions $f$ and $g$; see~\cite[Section 4.3]{Attouch2010} and~\cite[Lemma 2.3]{attouch2013convergence}.  Moreover, note that $\lambda_f=\infty$ and $g$ admits a Lipschitz continuous gradient with constant $L=1$; see~\cite[Exercise 1.24]{rockwets} and~\cite[Corollary 12.31]{BC}. The desired result then follows from Corollary~\ref{cor: semialgebraic convergence}.

(ii) Taking the fact that $0\in\partial F(x^*)$ into account and applying the subdifferential sum rule~\cite[Exercise 8.8]{rockwets}, one gets
\[x^*-\proj_C(x^*)\in-N_D(x^*) .\]
The constraint qualification~(\ref{CQ}) then implies that $x^*=\proj_C(x^*)$ and consequently $x^*\in C$. Moreover, the FRB scheme together with the closedness of $D$ ensures that $x^*\in D$, hence $x^*\in C\cap D$. Consequently, notice that the constraint qualification~(\ref{CQ}) amounts to
\[N_C(x^*)\cap\big(-N_D(x^*) \big)=\{0\}. \]
Then a direct application of~\cite[Theorem 5]{chen2020difference} guarantees that $F$ has KL exponent $\theta=\frac{1}{2}$ at $x^*$. The desired result follows immediately from Corollary~\ref{prop: KL exponent of merit function}.\endproof
\section{Numerical experiments}\label{num}
In this section, we apply the FRB splitting method to nonconvex feasibility problems.
Let $C=\{x\in\Rn:Ax=b\}$ for $A\in\R^{m\times n}$, $b\in\R^m$, $r=\ceil {m/5}$, $l=10^6$, and $D=\{x\in\Rn: \norm{x}_0\leq r, \norm{x}_\infty\leq l\}$. Similarly to~\cite[Section 5]{li2016douglas}, we consider the minimization problem~(\ref{optimization prblem}) with
\[f(x)=\delta_D\text{ and }g(x)=\frac{1}{2}\dist^2(x,C) .\]
Clearly $C$ is semialgebraic. As for $D$, first notice that $\Rn\times\{i\}$ for $i\in\R$ and
$\gph\norm{\cdot}_{0}$ are semialgebraic; see~\cite[Example 3]{Bolte2014}. Then
\begin{align*}
	\norm{x}_0\leq r&\Leftrightarrow \exists i\in\{0,\ldots,r\}\text{ such that }\norm{x}_0=i.\\
	&\Leftrightarrow \exists i\in\{0,\ldots,r\}\text{ such that }
x\in \proj_{\Rn}\big[\gph\norm{\cdot}_0\cap (\Rn\times\{i\})\big],
\end{align*}
which means that $\{x\in\Rn:\norm{x}_0\leq r \}$ is a finite union of intersections of semialgebraic sets, hence semialgebraic; see also~\cite[Formula 27(d)]{bauschke2014restricted}. On the other hand, one has $\norm{x}_\infty\leq l\Leftrightarrow\max_{1\leq i\leq n}(|x_i|-l)\leq0$, which means that the box $[-l,l]^n$ is semialgebraic. Altogether, the set $D$, which is intersection of semialgebraic sets, is semialgebraic. Hence, when specified to the problem above, FRB converges to a stationary point thanks to Example~\ref{ex: feasibility convergence}. %Moreover, $f$ is a proper lsc function bounded below, and $g$ is coercive and admits a Lipschitz continuous gradient with constant $L=1$; see, e.g.,~\cite[Corollary 12.31]{BC}.  Altogether, when specified to the problem above, FRB converges to a stationary point thanks to Example~\ref{ex: feasibility convergence}.

We find a projection onto $D$ by the formula given below, which is a consequence of~\cite[Proposition 3.1]{lu2013sparse} and was already observed by Li and Pong~\cite{li2016douglas}. We provide a proof for the sake of completeness.
\begin{proposition} Let $z=(z_1,\ldots,z_n)\in\Rn$. For every $i$, set $\tilde{z}_i^*=\proj_{[-l,l]}(z_i)$, $v_i^*=|z_i|^2-|\tilde{z}_i^*-z_i|^2$, and let $I^*\subseteq\{1,\ldots,n\}$ be the set of indices corresponding to the $r$ largest elements of $v_i^*,~i=1,\ldots,n$. Define $z^*\in\Rn$ by $z^*_i=\tilde{z}_i^*$ if $i\in I^*$ and $z^*_i=0$ otherwise. Then $z^*\in\proj_D(z)$.
\end{proposition}
\proof Apply \cite[Proposition 3.1]{lu2013sparse} with $\phi_i=|\cdot-z_i|^2$ and $\mathcal{X}_i=[-l,l]$.\endproof
%Moreover, invoking~\cite[Proposition 3.1]{lu2013sparse} provides

We shall benchmark FRB against
the Douglas-Rachford  method with fixed step-size~(DR)~\cite{li2016douglas} by Li and Pong, inertial Tseng's method~(iTseng)~\cite{boct2016inertial} by Bo{\c{t}} and Csetnek, and DR equipped with step-size heuristics~(DRh)~\cite{li2016douglas} by Li and Pong. These
splitting algorithms for nonconvex optimization problems are
known to converge globally to a stationary point of~(\ref{optimization prblem}) under appropriate assumptions on the concave KL property of  merit functions; see~\cite[Theorems 1--2, Remark 4, Corollary 1]{li2016douglas} and~\cite[Theorem 3.1]{boct2016inertial}. The convergence of DR and iTseng in our setting are already proved in~\cite[Proposition 2]{li2016douglas} and~\cite[Corollary 3.1]{boct2016inertial}, respectively.

We implement FRB with the following specified scheme for the problem of finding an $r$-sparse solution of a linear system $\{x\in\Rn:Ax=b\}$
\begin{align*}
&x_{k+1}\in\proj_D\left(x_k-\lambda A^\dagger A(2x_k-x_{k-1})+\lambda A^\dagger b \right)
%\argmin_{x\in D} \norm{x-x_k+\lambda A^\dagger A(2x_k-x_{k-1})+A^\dagger B } ,
\end{align*}
where the step-size $\lambda=0.9999\cdot\frac{1}{4}$ (recall Example~\ref{ex: feasibility convergence}). The inertial type Tseng's method studied in~\cite[Scheme (6)]{boct2016inertial} is applied with a step-size $\lambda^\prime= 0.1316$ given by~\cite[Lemma 3.3]{boct2016inertial} and a fixed inertial term $\alpha=\frac{1}{8}$:
\begin{align*}
%&p_{k+1}\in\argmin_{x\in D}\norm{x-x_k+\lambda^\prime A^\dagger(Ax_k-b)-\alpha(x_k-x_{k-1})},\\
&p_{k+1}\in\proj_D\left(x_k-\lambda^\prime A^\dagger(Ax_k-b)+\alpha(x_k-x_{k-1}) \right),\\
&x_{k+1}=p_k+\lambda^\prime A^\dagger A(x_k-p_k).
\end{align*}
As for DR and DRh, we employ the schemes specified by~\cite[Scheme (7)]{li2016douglas} and~\cite[Section 5]{li2016douglas} with the exact same step-sizes, respectively.  All algorithms are initialized at the origin, and we terminate FRB and iTseng when
\begin{equation*}
\frac{\max\big\{\norm{x_{k+1}-x_k },\norm{x_k-x_{k-1}} \big\} }{\max\big\{1,\norm{x_{k+1}},\norm{x_k},\norm{x_{k-1}}\big\}}<10^{-8}.
\end{equation*}
We adopt the termination criteria from~\cite[Section 5]{li2016douglas} for DR and DRh, where the same tolerance of $10^{-8}$ is applied. Similar to~\cite[Section 5]{li2016douglas}, our problem data is generated through creating random matrices $A\in\R^{m\times n}$ with entries following the standard Gaussian distribution. Then we generate a vector $\hat{x}\in\R^r$ randomly with the same distribution, project it onto the box $[-10^6,10^6]^r$, and create a sparse vector $\tilde{x}\in\Rn$ whose $r$ entries are chosen randomly to be the respective values of the projection of $\hat{x}$ onto $[-10^6,10^6]^r$. Finally, we set $b=A\tilde{x}$ to guarantee $C\cap D\neq\varnothing$.

Results of our experiments are presented in Table~\ref{num. results} below. For each problem of the size~$(m,n)$, we randomly generate 50 instances using the strategy described above, and report ceilings of the average number of iterations (iter), the minimal objective value at termination ($\text{fval}_{\text{min}}$), and the number of successes (succ). Here we say an experiment  is ``successful", if the objective function value at termination is less than $10^{-12}$, which means that the algorithm actually hits a global minimizer rather than just a stationary point. We observed the following:
\begin{itemize}
	\item[-] Among algorithms with fixed step-size (FRB, DR, and iTseng), FRB outperforms the others in terms of both the number of iterations and successes,  and it also has the smallest termination values. Moreover, it's worth noting that our simulation results align with Malitsky and Tam's observation that FRB converges faster than Tseng's method on a specific problem~\cite[Remark 2.8]{malitsky2020FRB}.
	\item[-] DRh has the most number of successes and the best precision at termination (see $\text{fval}_{\text{min}}$), but tend to be slower on ``easy" problems~(large $m$).
\end{itemize}
%Taking the above into account, we conclude that
Therefore, FRB is a competitive method for finding a sparse solution of a linear system, at least among the aforementioned algorithms with a fixed step-size\footnote{We also performed simulations with much larger problem size ($n=4000,5000,6000$), in which case FRB, DR, and iTseng tend to stuck at stationary points while DRh can still hit global minimizers. We believe that this is due to its heuristics; see also~\cite[Section 5]{li2016douglas}. }.

%Given the success of DRh on ``hard" problems, it would be tempting to see whether FRB equipped with appropriate heuristics has similar performance. We delay this investigation since the focus of this paper is to obtain global convergence of FRB to a stationary point, rather than exploring the best method for the nonconvex feasibility problem considered in this section.
\section{Conclusion and future work}\label{section: conclusion and future work}
We established convergence of the Malitsky-Tam FRB splitting algorithm in the nonconvex setting. Under the generalized concave KL property, we showed that FRB converges globally to a stationary point of~(\ref{optimization prblem}) and admit finite length property, which is even new in the convex setting. The sharpness of our results is demonstrated by virtue of the exact modulus associated with the generalized concave KL property. We also analyzed convergence rates of FRB when desingularizing functions associated with the generalized concave KL property have the \L ojasiewicz form. Numerical simulations suggest that FRB is a competitive method compared to DR and inertial Tseng's methods.

As for future work, it is tempting to analyze convergence rates using the exact modulus, as the usual KL convergence rates analysis assumes desingularizing functions of the \L ojasiewicz form, which could be an overestimation.  Another direction is to improve FRB through heuristics or other step-size strategies such as FISTA.

\section*{Acknowledgments}
XW and ZW were partially supported by NSERC Discovery Grants.

\begin{landscape}
\begin{table}
\caption{Comparing FRB to DR, iTseng, and DRh.\label{num. results}}
\begin{tabular}{|c c|c c c|c c c|c c c|c c c|}
\hline
Size   & &  FRB  &  &  & DR &  & & iTseng  &    && DRh  &   & \\
m    &n       &iter &$\text{fval}_{\text{min}}$&succ&iter  &$\text{fval}_{\text{min}}$&succ&iter  &$\text{fval}_{\text{min}}$&succ& iter  &$\text{fval}_{\text{min}}$& succ\\
\hline
300&600&411 &1.2756e-13                           &48    &476&1.6398e-13                         &43    &922&6.8569e-13                        &13      &436&1.0319e-30                         &50\\
300&700& 529& 1.4754e-13& 40&   601 &1.9085e-13& 36&1101& 7.9973e-13& 2& 444 &6.4897e-31 &50\\
300&800&665& 1.9931e-13& 29&    743 & 2.1789e-13& 22&1353& 1.2470e-12& 0& 448& 6.6553e-31& 50\\
300&900& 768&2.0614e-13&25&857&2.5572e-13&21&1537&1.2502e-12&0&452& 3.9219e-31&50\\
300& 1000&864&2.4851e-13 &16 & 963 &2.5456e-13& 11 & 1706 &1.3025e-12&0 & 457 &5.8254e-31&50\\
\hline
400&600 &238& 9.7199e-14 &50& 269& 1.1949e-13 &50&586& 5.7356e-13& 34&429 &4.4311e-30 &50\\
400&700&325& 1.0421e-13& 50&371&1.2968e-13&50&819 &5.5036e-13& 13&435 &1.9994e-30& 50\\
400&800&415 &1.7055e-13 &49&481 &1.9417e-13 &48&1045 &9.3313e-13 &2&439 &1.8757e-30& 50\\
400&900 &519 &2.1181e-13 &47  &    591 &2.2947e-13 &40   &   1183 &9.8694e-13 &1 &     442 &1.7107e-30 &50\\
 400 &    1000  &    609& 2.5329e-13 &40  &    688 &3.1439e-13 &32   &   1321 &1.4311e-12 &0   &   445&1.4946e-30 &50 \\
\hline
500&600&155& 9.0539e-14 &50&171 &1.0417e-13 &50&346& 5.2005e-13& 50&377& 1.6845e-29 &50\\
500&700&212 &1.2199e-13&50&239 &1.4644e-13& 50&509 &6.4501e-13& 28&421 &7.9475e-30 &50\\
500&800&273& 1.5619e-13 &50&310 &1.8283e-13& 50&688 &8.2500e-13& 10&432&5.0807e-30& 50\\
 500 & 900 &334 &1.7389e-13& 49&384& 2.1790e-13& 50 &863 &9.3320e-13 &3 &     433 &3.7216e-30 &50 \\
 500 & 1000& 414&2.1383e-13&50 & 474& 2.7150e-13 &49   &   1033 &1.1689e-12& 0   &   436 &2.8924e-30& 50\\
 \hline
\end{tabular}
\end{table}
\end{landscape}
\bibliographystyle{siam}
\bibliography{nonconvexFRB}

\end{document}